\documentclass{article}

\usepackage{amsmath}
\usepackage{amsthm}
\usepackage{alltt}
\usepackage[pdftex]{hyperref}

\numberwithin{equation}{section}

\newcommand{\res}[1]
{\textrm{Res}_{#1}}

\title{Some Variations of Two Combinatorial Identities}
\author{M.J. Kronenburg}
\date{}

\begin{document}

\maketitle

\begin{abstract}
Given two combinatorial identities proved earlier, a new set of variations of these
combinatorial identities is listed and proved with the integral representation method.
Some identities from literature are shown to be special cases of these new identities.
\end{abstract}


\section{Two Combinatorial Identities}

In an earlier paper \cite{K17}, the following two combinatorial identities were proved:
\begin{equation}\label{chugen}
 \sum_{k=0}^n \binom{a}{k}\binom{b}{n-k}\binom{k}{c}\binom{n-k}{d} 
  = \binom{a+b-c-d}{n-c-d}\binom{a}{c}\binom{b}{d}
\end{equation}
\begin{equation}\label{chu2gen}
 \sum_{k=0}^a \binom{a}{k}\binom{b}{m+k}\binom{k}{c}\binom{m+k}{d} 
  = \binom{a+b-c-d}{m+a-d}\binom{a}{c}\binom{b}{d}
\end{equation}
In this paper a list of identities that are variations of these
two identities is provided, and it is shown how they are proved
with the integral representation method.

\section{A Set of New Combinatorial Identities}

The following are variations of the first identity (\ref{chugen}):
\begin{equation}\label{eq1}
 \sum_{k=0}^n \binom{c+d-b}{k-p}\binom{b}{n-k}\binom{k}{c}\binom{n-k}{d} 
  = \binom{n-b}{n-d-p}\binom{p}{c+d+p-n}\binom{b}{d}
\end{equation}
\begin{equation}\label{eq2}
 \sum_{k=0}^a \binom{a}{k}\binom{c+d-a}{p+k}\binom{k}{c}\binom{n-k}{d} 
  = \binom{n+a+p-c-d}{a+p}\binom{n-a}{d-a-p}\binom{a}{c}
\end{equation}
\begin{equation}\label{eq3}
 \sum_{k=0}^a \binom{a}{k}\binom{c+d-a}{p-k}\binom{k}{c}\binom{n-k}{d} 
  = \binom{n-p}{c+d-p}\binom{n-a}{p-c}\binom{a}{c}
\end{equation}
\begin{equation}\label{eq4}
 \sum_{k=0}^n \binom{a}{k}\binom{b}{n-k}\binom{p+k}{a+b-d}\binom{n-k}{d} 
  = \binom{n+p-b}{n-d}\binom{p}{a+b-n}\binom{b}{d}
\end{equation}
\begin{equation}\label{eq5}
 \sum_{k=0}^n \binom{a}{k}\binom{b}{n-k}\binom{p-k}{a+b-d}\binom{n-k}{d} 
  = \binom{d+p-n}{a+b-n}\binom{p-a}{n-d}\binom{b}{d}
\end{equation}
\begin{equation}\label{eq6}
 \sum_{k=0}^n \binom{a}{k}\binom{b}{n-k}\binom{k}{c}\binom{p+k}{a+b-c} 
  = \binom{n+p-b}{n-c}\binom{p+c}{a+b-n}\binom{a}{c}
\end{equation}
\begin{equation}\label{eq7}
 \sum_{k=0}^n \binom{a}{k}\binom{b}{n-k}\binom{k}{c}\binom{p-k}{a+b-c} 
  = \binom{p-n}{a+b-n}\binom{p-a}{n-c}\binom{a}{c}
\end{equation}
The following are variations of the second identity (\ref{chu2gen}):
\begin{equation}\label{eq8}
 \sum_{k=0}^{a+p} \binom{a}{k-p}\binom{b}{m+k}\binom{k}{c}\binom{m+k}{m} 
  = \binom{a+b-c-m}{a+p-c}\binom{b-m}{c}\binom{b}{m}
\end{equation}
\begin{equation}\label{eq9}
 \sum_{k=0}^a \binom{a}{k}\binom{c+d-a}{p+k}\binom{k}{c}\binom{m+k}{d} 
  = \binom{m-p}{d-a-p}\binom{m+c}{a+p}\binom{a}{c}
\end{equation}
\begin{equation}\label{eq10}
 \sum_{k=0}^a \binom{a}{k}\binom{c+d-a}{p-k}\binom{k}{c}\binom{m+k}{d} 
  = \binom{m+a+p-c-d}{p-c}\binom{m+c}{c+d-p}\binom{a}{c}
\end{equation}
\begin{equation}\label{eq11}
 \sum_{k=0}^a \binom{a}{k}\binom{b}{m+k}\binom{p+k}{a+b-d}\binom{m+k}{d} 
  = \binom{d+p-m}{b-m}\binom{p}{m+a-d}\binom{b}{d}
\end{equation}
\begin{equation}\label{eq12}
 \sum_{k=0}^a \binom{a}{k}\binom{b}{m+k}\binom{p-k}{a+b-d}\binom{m+k}{d} 
  = \binom{m+p-b}{m+a-d}\binom{p-a}{b-m}\binom{b}{d}
\end{equation}
\begin{equation}\label{eq13}
 \sum_{k=0}^a \binom{a}{k}\binom{b}{m+k}\binom{k}{c}\binom{p+k}{a+b-c} 
  = \binom{p-m}{b-c-m}\binom{p+c}{m+a}\binom{a}{c}
\end{equation}
\begin{equation}\label{eq14}
 \sum_{k=0}^a \binom{a}{k}\binom{b}{m+k}\binom{k}{c}\binom{p-k}{a+b-c} 
  = \binom{m+p-b}{m+a}\binom{p-a}{b-m-c}\binom{a}{c}
\end{equation}

\section{Proof of the New Combinatorial Identities}
The new combinatorial identities listed above can all be proved with the
integral representation method \cite{CB84,E84}.
Two of these identities are proved below, and the other identities
have similar proofs.\\
For identity (\ref{eq1}), applying the trinomial revision identity \cite{GKP94,K97,K15},
the left side of the identity simplifies to:
\begin{equation}
 \sum_{k=0}^n \binom{a}{k-p}\binom{k}{c}\binom{b-d}{n-d-k} 
\end{equation}
Using the integral representation method \cite{CB84,E84}, this becomes:
\begin{equation}
\begin{split}
 & \sum_{k=0}^{\infty} \res{x}\frac{(1+x)^a}{x^{k-p+1}}\res{y}\frac{(1+y)^k}{y^{c+1}}\res{z}\frac{(1+z)^{b-d}}{z^{n-d-k+1}} \\
 ={} & \res{x}\res{y}\res{z} \frac{(1+x)^a(1+z)^{b-d}}{x^{-p+1}y^{c+1}z^{n-d-1}} \sum_{k=0}^{\infty}\left(\frac{(1+y)z}{x}\right)^k \\
 ={} & \res{y}\res{z}\res{x} \frac{(1+x)^a(1+z)^{b-d}x^{p}}{y^{c+1}z^{n-d+1}(x-(1+y)z)} \\
 ={} & \res{y}\res{z} \frac{(1+(1+y)z)^a(1+y)^p(1+z)^{b-d}}{y^{c+1}z^{n-d-p+1}} \\
\end{split}
\end{equation}
Now the following is used:
\begin{equation}
 (1+(1+y)z)^a = ((1+z)+yz)^a = \sum_{k=0}^a \binom{a}{k}(1+z)^k (yz)^{a-k}
\end{equation}
Collecting the residues, the following expression results:
\begin{equation}
\begin{split}
 & \sum_{k=0}^a \binom{a}{k}\binom{p}{c-a+k}\binom{b-d+k}{n-d-p-a+k} \\
 ={} & \sum_{k=0}^a \binom{a}{k}\binom{p}{c-a+k}\binom{b-d+k}{a+b+p-n} \\
\end{split}
\end{equation}
This expression is recognized as (\ref{chu2gen}) with $c=0$ provided that $m=c-a=b-d$,
or equivalently $a=c+d-b$. The parameters of (\ref{chu2gen}) thus become
$a\leftarrow c+d-b$, $b\leftarrow p$, $c\leftarrow 0$, $d\leftarrow c+d+p-n$ and $m\leftarrow b-d$,
and substituting these parameters in the right side of (\ref{chu2gen}) gives the result.\\
For identity (\ref{eq2}), applying the trinomial revision identity \cite{GKP94,K97,K15},
the left side of the identity simplifies to:
\begin{equation}
 \sum_{k=0}^n \binom{a-c}{k-c}\binom{b}{p+k}\binom{n-k}{d} 
\end{equation}
Using the integral representation method \cite{CB84,E84}, this becomes:
\begin{equation}
\begin{split}
 & \sum_{k=0}^{\infty} \res{x}\frac{(1+x)^{a-c}}{x^{k-c+1}}\res{y}\frac{(1+y)^b}{y^{p+k+1}}\res{z}\frac{(1+z)^{n-k}}{z^{d+1}} \\
 ={} & \res{x}\res{y}\res{z} \frac{(1+x)^{a-c}(1+y)^b(1+z)^n}{x^{-c+1}y^{p+1}z^{d+1}} \sum_{k=0}^{\infty}\left(\frac{1}{xy(1+z)}\right)^k \\
 ={} & \res{y}\res{z}\res{x} \frac{(1+x)^{a-c}(1+y)^b(1+z)^n x^{c}}{y^{p+1}z^{d+1}(x-1/(y(1+z)))} \\
 ={} & \res{y}\res{z} \frac{(1+1/(y(1+z)))^{a-c}(1+y)^b(1+z)^{n-c}}{y^{c+p+1}z^{d+1}} \\
 ={} & \res{y}\res{z} \frac{(1+y(1+z))^{a-c}(1+y)^b(1+z)^{n-a}}{y^{a+p+1}z^{d+1}} \\
\end{split}
\end{equation}
Now the following is used:
\begin{equation}
 (1+y(1+z))^{a-c} = ((1+y)+yz)^{a-c} = \sum_{k=0}^{a-c} \binom{a-c}{k}(1+y)^k (yz)^{a-c-k}
\end{equation}
Collecting the residues, the following expression results:
\begin{equation}
\begin{split}
 & \sum_{k=0}^{a-c} \binom{a-c}{k}\binom{b+k}{c+p+k}\binom{n-a}{c+d-a+k} \\
 ={} & \sum_{k=0}^{a-c} \binom{a-c}{k}\binom{b+k}{b-c-p}\binom{n-a}{c+d-a+k} \\
\end{split}
\end{equation}
This expression is recognized as (\ref{chu2gen}) with $c=0$ provided that $m=b=c+d-a$.
The parameters of (\ref{chu2gen}) thus become
$a\leftarrow a-c$, $b\leftarrow n-a$, $c\leftarrow 0$, $d\leftarrow d-p-a$ and $m\leftarrow c+d-a$,
and substituting these parameters in the right side of (\ref{chu2gen}) gives the result.\\

\section{Some Identities from Literature}
Some identities from literature are special cases of the identities above.\\
T.S. Nanjundiah gave the following two identities \cite{B70,G72,N58,Q16}:
\begin{equation}\label{nanj1}
 \sum_{k=0}^a \binom{a}{k}\binom{b}{k}\binom{p+k}{a+b} = \binom{p}{a}\binom{p}{b}
\end{equation}
\begin{equation}\label{nanj2}
 \sum_{k=0}^n \binom{m-x+y}{k}\binom{n+x-y}{n-k}\binom{x+k}{m+n} = \binom{x}{m}\binom{y}{n}
\end{equation}
The first equation is equation (\ref{eq11}) with $m=d=0$,
and the second equation is equation (\ref{eq4}) with $a=m-x+y$, $b=n+x-y$, $d=0$ and $p=x$.\\
M.T.L. Bizley gave the following two identities \cite{B70,G72}:
\begin{equation}\label{biz1}
 \sum_{k=0}^a \binom{a}{k}\binom{b}{k-d}\binom{p+k}{a+b} = \binom{p}{a-d}\binom{p+d}{b+d}
\end{equation}
\begin{equation}\label{biz2}
 \sum_{k=0}^n \binom{a}{k}\binom{b}{n-k}\binom{p+k}{a+b} = \binom{p}{a+b-n}\binom{p-b+n}{n}
\end{equation}
The first equation is equation (\ref{eq9}) with $c=0$, $d=a+b$, $m=p$ and $p=-d$,
and the second equation is equation (\ref{eq4}) with $d=0$.\\
H.W. Gould gave the following identity \cite{G72}:
\begin{equation}\label{sur1}
 \sum_{k=0}^a \binom{a}{k}\binom{b}{k}\binom{a+b+x+k}{a+b} = \binom{a+b+x}{a}\binom{a+b+x}{b}
\end{equation}
This equation is equation (\ref{eq11}) with $m=d=0$ and $p=a+b+x$.\\
J. Sur\'anyi gave the following identity \cite{G72,Q16,S85}:
\begin{equation}\label{sur1}
 \sum_{k=0}^a \binom{a}{k}\binom{b}{k}\binom{a+b+x-k}{a+b} = \binom{x+a}{a}\binom{x+b}{b}
\end{equation}
This equation is equation (\ref{eq12}) with $m=d=0$ and $p=a+b+x$.\\
L. Tak\'acs gave the following identity \cite{S85,T73}:
\begin{equation}\label{tak1}
 \sum_{k=0}^n \binom{a}{k}\binom{m-a}{n-k}\binom{p+k}{m} = \binom{p}{m-n}\binom{n+a+p-m}{n}
\end{equation}
This equation is equation (\ref{eq4}) with $b=m-a$ and $d=0$.\\
J. Riordan gave the following identity \cite{G72}:
\begin{equation}\label{tak1}
 \sum_{k=0}^n \binom{n}{k}\binom{m}{n-k}\binom{x+n-k}{n+m} = \binom{x}{m}\binom{x}{n}
\end{equation}
This equation is equation (\ref{eq5}) with $a=n$, $b=m$, $d=0$ and $p=x+n$.\\
R.P. Stanley gave the following two identities \cite{G72,G72a,S71}:
\begin{equation}\label{stan1}
 \sum_{k=0}^a \binom{p+q+k}{k}\binom{p}{a-k}\binom{q}{b-k} = \binom{p+b}{a}\binom{q+a}{b}
\end{equation}
\begin{equation}\label{stan2}
 \sum_{k=0}^a (-1)^k \binom{p+q+1}{k}\binom{p+a-k}{p}\binom{q+b-k}{q} = \binom{p+a-b}{a}\binom{q+b-a}{b}
\end{equation}
These two identities can be shown to be equivalent to (\ref{eq5}) with $d=0$
using the following two symmetry identities \cite{GKP94,K97,K15}:
\begin{equation}
 \binom{n}{k} = (-1)^k \binom{-n+k-1}{k} = (-1)^{n-k}\binom{-k-1}{n-k}
\end{equation}
In the first identity, applying the first symmetry identity
to the first and third binomial coefficients on the left and to the second
binomial coefficient on the right, and then replacing $q$ by $-q-1$,
results in the following identity:
\begin{equation}\label{stanres}
 \sum_{k=0}^a \binom{q-p}{k}\binom{p}{a-k}\binom{q+b-k}{q} = \binom{p+b}{a}\binom{q+b-a}{b}
\end{equation}
When $q<p$, interchanging $p$ with $q$ and $a$ with $b$, which does not change Stanley's identities, makes $q>p$.
This equation is equation (\ref{eq5}) with $a\leftarrow q-p$, $b\leftarrow p$, $d\leftarrow 0$, $n\leftarrow a$ and $p\leftarrow q+b$.\\
In the second identity, applying the second symmetry identity
to the second binomial coefficient on the left and the first symmetry identity
to the first binomial coefficient on the right, and then replacing $p$ with $-p-1$,
results in the same identity.

\pdfbookmark[0]{References}{}


\begin{thebibliography}{99}
\bibitem{B70}
  M.T.L. Bizley,
  A Generalization of Nanjundiah's Identity,
  \textit{Amer. Math. Monthly} 77~(1970)~863-865.
\bibitem{CB84}
  R.V. Churchill, J.W. Brown,
  \textit{Complex Variables and Applications},
  McGraw-Hill, 1984.
\bibitem{E84}
  G.P. Egorychev,
  \textit{Integral Representation and the Computation of Combinatorial Sums},
  Translations of Mathematical Monographs, 59, Amer. Math. Soc., 1984.
\bibitem{G72}
  H.W. Gould,
  \textit{Combinatorial Identities}, rev. ed.,
  Morgantown, 1972.
\bibitem{G72a}
  H.W. Gould,
  A New Symmetrical Combinatorial Identity,
  \textit{J. Combin. Theory Ser. A} 13~(1972)~278-286.
\bibitem{GKP94}
  R.L. Graham, D.E. Knuth, O. Patashnik,
  \textit{Concrete Mathematics, A Foundation for Computer Science}, 2nd ed.,
  Addison-Wesley, 1994.
\bibitem{K97}
  D.E. Knuth,
  \textit{The Art of Computer Programming, Volume 1: Fundamental Algorithms}, 3rd ed.,
  Addison-Wesley, 1997.
\bibitem{K15}
  M.J. Kronenburg,
  The Binomial Coefficient for Negative Arguments,
  \href{http://arxiv.org/abs/1105.3689}{{\tt arXiv:1105.3689}}{\tt~[math.CO]}
\bibitem{K17}
  M.J. Kronenburg,
  A Generalization of the Chu-Vandermonde Convolution and some Harmonic Number Identities,
  \href{http://arxiv.org/abs/1701.02768}{{\tt arXiv:1701.02768}}{\tt~[math.CO]}
\bibitem{N58}
  T.S. Nanjundiah,
  Remark on a note of P. Tur\'an,
  \textit{Amer. Math. Monthly} 65~(1958)~354.
\bibitem{Q16}
  J. Quaintance, H.W. Gould,
  \textit{Combinatorial Identities for Stirling Numbers},
  World Scientific, 2016.
\bibitem{S85}
  L.A. Sz\'ekely,
  Common Origin of Cubic Binomial Identities;
  A Generalization of Sur\'anyi's Proof on
  Le Jen Shoo's Formula,
  \textit{J. Combin. Theory Ser. A} 40~(1985)~171-174.
\bibitem{S71}
  R.P. Stanley,
  \textit{Ordered Structures and Partitions},
  Ph. D. thesis, Harvard University, 1971.
\bibitem{T73}
  L. Tak\'acs,
  On an identity of Shih-Chieh Chu,
  \textit{Acta Sci. Math. (Szeged)} 34~(1973)~383-391.
\end{thebibliography}
\end{document}